\documentclass[reqno]{amsart}
\usepackage{latexsym}
\usepackage{amssymb,amsthm,amsmath}
\usepackage[dvips]{graphicx}
\usepackage{changebar}
\allowdisplaybreaks 

\usepackage{kantlipsum} 

\calclayout

\theoremstyle{plain}
\newtheorem{theorem}{Theorem}

\theoremstyle{definition}
\newtheorem{definition}[theorem]{Definition}
\theoremstyle{remark}

\def\pa{\partial}
\def\cal{\mathcal}
\let\mib=\boldsymbol

\def\R{{\mathbb R}}
\def\N{{\mathbb N}}

\def\eps{\varepsilon}

\def\malpha{{\mib \alpha}}
 \def\mbeta{{\mib \beta}}

\def\mx{{\bf x}}
\def\mz{{\bf z}}
\def\mxi{{\mib \xi}}
\def\my{{\bf y}}
\def\mz{{\bf z}}
\def\mw{{\bf w}}

\begin{document}

\title[]{The structure of ${\cal A}$-free measures with uniformly singular part}
\author{D. Mitrovic}\address{Darko Mitrovic,
Faculty of Mathematics, University of
Montenegro 81000 Podgorica, Montenegro}\email{darkom@ac.me}


\begin{abstract}
We prove that a singular part $\mu_s$ of a measure $\mu$ satisfying ${\cal A}\mu =0$ for a linear partial differential operator ${\cal A}$ defined on $\R^d$ has the range in the intersection of kernels of the principal symbol of ${\cal A}$ if the singular part is singular with respect to all the variables (uniformly singular) i.e. it is such that for $\mu_s$-almost every $\mx\in \R^d$ there exist positive functions $\alpha(\eps), \beta(\eps)$, $\eps \in \R$, satisfying $\frac{\alpha(\eps)}{\eps}\to 0$, $ \frac{\eps}{\beta(\eps)}\to 0$ and a set $E_\eps\subset B(\mx,\alpha(\eps))$ such that $\lim\limits_{\eps\to 0}\frac{\mu_s(B(\mx,\beta(\eps))\backslash E_\eps)}{|\mu_s|(E_\eps)}=0$.

\end{abstract}
\maketitle

\section{Introduction}

In the paper, we consider a finite Radon measure $\mu=(\mu_1,\dots, \mu_m)$ defined on $\R^d$ satisfying the system of partial differential equation
\begin{equation}
\label{main-eq}
{\cal A} \mu=\sum\limits_{\malpha \in I} \pa^\malpha (A_\malpha \mu )=0,
\end{equation} where $I=I_1\times I_2\times\dots \times I_n\subset\{\malpha=(\alpha_1,\dots,\alpha_d):\, \alpha_s\in \N\cup \{0\}, \, s=1,\dots,d \}^n$ is a set of multi-indexes, $\pa^\malpha=\pa_{x_1}^{\alpha_1} \pa_{x_2}^{\alpha_2}\dots \pa_{x_d}^{\alpha_d}$, and 
$A_\malpha:\R^d \to M^{n\times m}$ are smooth mappings from $\R^d$ into the space of real ${n\times m}$ matrices. Written coordinate-wise, we actually have the following system of equations 
\begin{equation}
\label{j}
{\cal A}_j \mu= \sum\limits_{\malpha \in I_j} \sum\limits_{k=1}^m \pa^\malpha (a^{\malpha}_{jk}\mu_k)=0, \ \ j=1,\dots,n,
\end{equation} where $I_j\subset \{\malpha=(\alpha_1,\dots,\alpha_d):\, \alpha_s\in \N\cup \{0\}, \, s=1,\dots,d \}$. Denote by $A_j$, $j=1,\dots,n$, the principal symbol of the operator ${\cal A}_j$ given by
\begin{equation}
\label{glavni}
 A_j(\mx, \mxi)=\sum\limits_{\malpha \in I_j'}  \sum\limits_{k=1}^m a^\malpha_{jk}(\mx) (2\pi
i\mxi)^{\malpha}, \ \ I_j' \subset I_j.
\end{equation} The sum given above is taken over all terms from \eqref{j}
whose order of derivative $\malpha$ is not dominated by any other
multi-index from $I_j$. As usual, $\mxi^\malpha=\xi_1^{\alpha_1}\dots\xi_d^{\alpha_d}$ for $\malpha=(\alpha_1,\dots,\alpha_d)$, and $|\malpha|=\alpha_1+\dots+\alpha_d$.

For instance, for the (scalar) operator ${\cal A}=\pa_{x_1}+\pa_{x_2}+\pa^2_{x_2}$, we have $I=I_1=\{ (1,0),(0,1), (0,2) \}$ and $I'=I'_1=\{(1,0),(0,2) \}$.

We are interested in the range of the Radon-Nikodym derivative $f(\mx)=\frac{d \mu_s}{d |\mu_s|}(\mx)$ of the singular part $\mu_s$ of the measure $\mu=\mu_a+\mu_s$ where the latter is the Lebesgue decomposition of the measure $\mu$. The problem is initiated in \cite{AdG} where it is conjectured that for the $k$-th order operator ${\cal A}$, the function $f$ must take values in the wave cone $\Lambda_{{\cal A}}=\cup_{|\mxi|=1} {\rm Ker}{\bf A}^k(\mxi)$ where ${\bf A}^k(\mxi)$ is the sum of all symbols of order $k$ (see \cite{DpR} for details). The problem is resolved in \cite{DpR} where one can find thorough information on this issue concerning history and applications (in particular in the calculus of variations and geometric measure theory). 

Here, we shall improve the result from \cite{DpR} by describing behavior of $\mu_s$ on more general manifolds. Moreover, we prove a stronger statement in the sense that the support of $f$ is actually not in the union of kernels but in their intersection (of course, if we assume that $\mu_s$ is uniformly singular; see Definition \ref{unif-s}).

To this end, for every of the principal symbols $A_j$, $j=1,\dots,n$, we  assume that there exists a multi-index
$\mbeta^j=(\beta^j_1, \cdots,
 \beta^j_d)\in I_j' \subset \N^d_+$ such that for any positive $\lambda\in \R$  the
following homogeneity assumption holds
\begin{equation}
\label{gen_homog}
A_j(\mx,\lambda^{\beta^j_1}\xi_1,\dots,\lambda^{\beta^j_d}\xi_d)=
\lambda A_j(\mx, \mxi),
\end{equation}
 implying that
 \begin{equation}
 \label{sum1}
 \sum\limits_{k=1}^d {\alpha_{k}}{\beta^j_k}=1 \ \ \text{for every $\malpha=(\alpha_1,\dots,\alpha_d) \in I'_j$}.
 \end{equation} We then introduce the homogeneity manifolds:
 $$
 P_j=\{ \mxi \in \R^d: \, |\xi_1|^{1/\beta^j_1}+\dots+|\xi_d|^{1/\beta^j_d}=1 \}
 $$ and the corresponding projections
 $$
 \pi_j(\mxi)\!=\!\left( \frac{\xi_1}{\left(|\xi_1|^{1/\beta^j_1}+\dots+|\xi_d|^{1/\beta^j_d}\right)^{\beta^j_1}},\dots, \frac{\xi_d}{\left(|\xi_1|^{1/\beta^j_1}+\dots+|\xi_d|^{1/\beta^j_d}\right)^{\beta^j_d}}\right)\!\!, \ \ \mxi \in \R^d.
 $$ In the case of the operator ${\cal A}=\pa_{x_1}+\pa_{x_2}+\pa^2_{x_2}$, the corresponding symbol is $A(\xi_1,\xi_2)=-i(\xi_1+\xi_2)+\xi_2^2$, the principal symbol is $A'(\xi_1,\xi_2)=-i\xi_1+\xi_2^2$, and $\mbeta=(1,1/2)$.
 
Finally, we need a condition on the singular part of the measure $\mu$ which we call the uniform singularity condition. Roughly speaking, we require that $\mu_s$ is singular with respect to every of the variables. For instance, such a condition is not fulfilled by the measure $\delta(x_1)dx_2$ since it is not singular with respect to $x_2$. 

\begin{definition} 
\label{unif-s}
We say that the measure $\mu_s$ is uniformly singular if for $\mu_s$-almost every $\mx\in \R^d$ there exist real positive functions $\alpha(\eps), \beta(\eps)$, $\eps \in \R$, satisfying $\frac{\alpha(\eps)}{\eps}\to 0$, $\frac{\eps}{\beta(\eps)}\to 0$, and a family of sets $E_\eps\subset B(\mx,\alpha(\eps))$ such that $\lim\limits_{\eps\to 0}\frac{\mu_s(B(\mx,\beta(\eps))\backslash E_\eps)}{|\mu_s|(E_\eps)}=0$.
\end{definition} For instance, it is clear that the measure $\delta(x_1)\delta(x_2)$ satisfies the latter condition.

We have the following theorem.
 
 \begin{theorem}
 \label{thm:main}
 Let $\mu$ be a solution to \eqref{main-eq} and let $\mu=h(\mx)d\mx+\mu^s$ be the Lebesgue decomposition of the measure $\mu$. Assume that the measure $\mu_s$ is uniformly singular in the sense of Definition \ref{unif-s}. Denote by $f=(f_1,\dots,f_m) \in L_{|\mu^s|}^1(\R^d;\R^m)$ the Radon-Nykodim derivatives of $\mu^s$ with respect to its total variation measure $|\mu^s|$: $f=\frac{d \mu^s}{ d |\mu^s|}$. It holds for 
 $|\mu^s|$-almost every $\mx\in \R^d$:
 \begin{equation}
 \label{kernel}
 \sum\limits_{\malpha \in I_j'} \sum\limits_{k=1}^m a^\malpha_{jk}(\mx) (2\pi \, i\mxi)^{\malpha} f_k(\mx)=0, \ \ {\cal H}^{d-1}-{\rm a.e.} \ \ \mxi \in  P_j, \ \ j=1,\dots,m.
 \end{equation}
  \end{theorem}
Remark that each equation of system \eqref{j} defines different homogeneity manifold. If all the manifolds $P_j$, $j=1,\dots,n$, would be the same, say $P$, and we have the same set of dominating multi-indices $I'=I'_j$, $j=1,\dots,n$, for then, denoting $\pi=\pi_j$, $j=1,\dots,n$, we could rewrite \eqref{kernel} in the form
$$
\sum\limits_{\malpha \in I'} (2 \pi \, i \mxi)^\malpha A_\malpha(\mx) f(\mx)=0 \ \ {\cal H}^{d-1}-{\rm a.e.} \ \ \mxi \in  P_j, \ \ j=1,\dots,m,
$$for appropriate matrices $A_\malpha$, $\malpha \in I'$. From here, we see that in the latter case, the statement of the theorem is actually
\begin{equation}
\label{comp-RD}
f(\mx)\in \cap_{\mxi \in P} {\rm Ker}\left(\sum\limits_{\malpha \in I'} (2 \pi \, i \mxi)^\malpha A_\malpha(\mx) \right) \ \ {\rm for} \ \ \mu_s \ \ {a.e.} \ \ \mx\in \R^d.
\end{equation} We remark that in \cite{DpR}, a constant coefficients operator ${\cal A}$ of order $k\in \N$ is considered and it was proved that \eqref{comp-RD} holds with the union $\Cup$ (instead of $\cap$) for $|\malpha|=k$ (instead of $\malpha \in I'$) with $|\mxi|=1$ (instead of $\mxi \in P_j$, $j=1,\dots, m$).
  
We will dedicate the last section to the proof of the theorem. In the next section, we shall prove it in the case of first order constant coefficients operator and the scalar measure which captures all the elements of the general situation.  The proof is based on the blow up method \cite{Vass} and appropriate usage of Fourier multiplier operators (as in deriving appropriate defect functionals \cite{AM, LM2}).

Let us recall that the Fourier multiplier operator $T_\psi$ with the symbol $\psi$ is defined vie the Fourier and inverse Fourier transform 
$$
T_\psi u (\mx)={\cal F}^{-1}(\psi(\mxi){\cal F}(u))(\mx), \ \ u\in L^2(\R^d),
$$where the Fourier and the inverse Fourier transforms are given by 
$$
{\cal F}(u)(\mxi)=\hat{u}(\mxi)=\int e^{-2\pi i \mx\cdot \mxi}u(\mx)d\mx, \ \ {\cal F}^{-1}(u)(\mx)=\check{u}(\mx)=\int e^{2\pi i \mx\cdot \mxi}u(\mxi)d\mxi.
$$ For properties of the Fourier multiplier operators one can consult \cite{G}.


\section{Proof of Theorem \ref{thm:main} in the case of first order constant coefficients operator and the scalar measure}

Here, we shall prove Theorem \ref{thm:main} when the scalar finite Radon measure $\mu\in {\cal M}(\R^d)$ satisfies the equation
\begin{equation}
\label{scalar}
\sum\limits_{l=1}^d \pa_{x_l}( a_l \mu) +a_0 \mu=0,
\end{equation}where $(a_0,a_1,\dots,a_d)$ is a constant vector. The proof is essentially the same for the general operator of the form given in \eqref{main-eq}, but the proof is a bit less technical for \eqref{scalar}. The proof in full generality is given in the next section.

Before we start, let $f$ be the Radon-Nykodim derivative of $\mu^s$ with respect to $|\mu^s|$ (we disregard the fact that $f$ can take only values $\pm 1$):
$$
d\mu^s(\my)=f(\my) d|\mu^s|(\my).
$$

We fix a convolution kernel $\rho:\R^d \to \R$ which is a smooth, compactly supported function of total mass one and convolve \eqref{scalar} by $\rho_\eps(\mx)=\frac{1}{\eps^d}\rho(\frac{\mx}{\eps})$. Then, we take an arbitrary $\varphi\in C^1_c(\R^d)$ and test the convolved equation on such a function. We get (below we denote $\mu^\eps=\mu\star \rho_\eps$ and $\langle \varphi(\my),\mu(\my) \rangle = \int_{\R^d} \varphi(\my)d\mu(\my)$)

\begin{equation}
\label{i-1}
\int_{\R^d} \langle \frac{1}{\eps^d}\rho\big(\frac{\mx\!-\!\my}{\eps}\big),\mu(\my) \rangle \sum\limits_{l=1}^d a_l \pa_{x_l}\varphi(\mx)d\mx+ a_0\int_{\R^d} \varphi(\mx) d\mu^\eps(\mx)=0.
\end{equation} We now fix $\mz \in \R^d$ and take $\varphi(\frac{\mx-\mz}{\eps})$ instead of $\varphi$ in \eqref{i-1}. We get:

\begin{equation}
\label{i-2}
\int_{\R^d} \langle \frac{1}{\eps^d}\rho\big(\frac{\mx\!-\!\my}{\eps}\big),\mu(\my) \rangle \sum\limits_{l=1}^d a_l \frac{1}{\eps} \pa_{w_l}\varphi(\mw)\Big|_{\mw=\frac{\mx-\mz}{\eps}}d\mx+ a_0\int_{\R^d} \varphi\big(\frac{\mx\!-\!\mz}{\eps}\big) d\mu^\eps(\mx)=0.
\end{equation}We now introduce in the first integral the change of variables $\mx=\mz+\eps \mw$ and multiply entire expression by $\eps$. We get:
\begin{equation}
\label{i-3}
\int_{\R^d} \langle \rho\big(\frac{\mz\!-\!\my}{\eps}\!+\!\mw\big),\mu(\my) \rangle \sum\limits_{l=1}^d a_l \pa_{w_l}\varphi(\mw)d\mw+ \eps a_0\int_{\R^d} \varphi\big(\frac{\mx\!-\!\mz}{\eps}\big) d\mu^\eps(\mx)=0.
\end{equation} Now, we shall use the uniform singularity conditions. We rewrite \eqref{i-3} in the form

\begin{equation}
\label{i-4}
\left(\int_{E_\eps}+\int_{B(\mz,\beta(\eps))\backslash E_{\eps}}+ \int_{\R^d\backslash B(\mz,\beta(\eps))} \right) \langle \rho\big(\frac{\mz\!-\!\my}{\eps}\!+\!\mw\big),\mu(\my) \rangle \sum\limits_{l=1}^d a_l \pa_{w_l}\varphi(\mw)d\mw+ \eps a_0\int_{\R^d} \varphi\big(\frac{\mx\!-\!\mz}{\eps}\big) d\mu^\eps(\mx)=0
\end{equation} for the set $E_\eps$ and the function $\beta(\eps)$ corresponding to the point $\mz$ in Definition \ref{unif-s}. Now, according to the assumptions for the uniformly singular measures and the fact that $\rho$ is compactly supported:
$$
E_\eps\subset B(\mz,\alpha(\eps)), \; \frac{\alpha(\eps)}{\eps}\to 0; \ \ \frac{\mu_s(B(\mz,\beta(\eps))\backslash E_{\eps})}{|\mu_s|(E_\eps)}\to 0; \ \ \rho(\frac{\mz-\my}{\eps}+\mw)\equiv 0, \; \mz\notin B(\mz,\beta(\eps)),
$$we get after dividing \eqref{i-4} by $|\mu_s|(E_\eps)$ and letting $\eps\to 0$ in \eqref{i-4}

\begin{equation}
\label{i-6}
\int_{\R^{d}} \rho(\mw)   \sum\limits_{l=1}^d a_l \pa_{w_l}\varphi(\mw) f(\mz) d\mw  =0.
\end{equation} Now, by simple density arguments, we see that we can take $\varphi\in C^1_0(\R^d)$. We choose
$$
\varphi(\mw)=\overline{T_{_{\!\!\frac{\psi(\mxi/|\mxi|)}{|\mxi|}}}\rho(\mw)},
$$where $T_{\frac{\psi(\mxi/|\mxi|)}{|\mxi|}}$ is the Fourier multiplier operator with the symbol ${\frac{\psi(\mxi/|\mxi|)}{|\mxi|}}$. We get from \eqref{i-6} after applying the Plancherel theorem with respect to $\mw$:

\begin{equation}
\label{i-7}
f(\mz) \int_{\R^{d}}    \sum\limits_{l=1}^d a_l \frac{2\pi \, i \xi_l}{|\mxi|}  \psi(\mxi/|\mxi|) |\hat{\rho}(\mxi)|^2   d\mxi =0.
\end{equation} From here, due to arbitrariness of $\psi$ and $\mz$, we conclude that for $|\mu^s|$-a.e. $\mz\in \R^d$,
$$
\sum\limits_{l=1}^d a_l \xi_l f(\mz)=0 \ \ \text{${\cal H}^{d-1}$-a.e. $\mxi=(\xi_1,\dots,\xi_d)$, $|\mxi|=1$}
$$what we wanted to prove.

\section{Proof of Theorem \ref{thm:main}; general case}

We follow the steps from the previous section and we address a reader there for clarifications.

We start by fixing $j$ in \eqref{j} and the convolution kernel $\rho:\R\to \R$ which is smooth, compactly supported with total mass one. We then denote ($\mx=(x_1,\dots,x_d)$ and $\mw=(w_1,\dots,w_d)$ below)
\begin{align*}
\rho_{j,\eps}(\mx)=\frac{1}{\eps^{\beta_1^j+\dots+\beta_d^j}}\prod\limits^d_{k=1} \rho(\frac{x_j}{\eps^{\beta_k^j}})\qquad {\rm and} \qquad
{\mib \rho}(\mw)=\prod\limits^d_{k=1}\rho(w_k),
\end{align*} and convolve \eqref{j} by $\rho_{j,\eps}$. We have for $(a^{\malpha}_{jl}\mu^s)^\eps=(a^{\malpha}_{jl}\mu_l)\star \rho_{j,\eps}(\mx)$:
\begin{equation}
\label{j-eps}
\sum\limits_{\malpha \in I_j}\sum\limits_{l=1}^m  \pa^\malpha (a^{\malpha}_{jl}\mu_l)^\eps=0.
\end{equation} We then apply a test function $\varphi \in C^\infty_c(\R^d)$ on \eqref{j-eps} to get
\begin{equation}
\label{j-1}
\sum\limits_{\malpha \in I_j}(-1)^{|\malpha|}\sum\limits_{l=1}^m \int_{\R^d} \langle \rho_{j,\eps}(\mx-\my), a^{\malpha}_{jl}(\my) \mu_l(\my) \rangle \pa^\malpha \varphi(\mx) d\mx=0.
\end{equation} Now, we fix $\mz\in \R^d$ and take $\varphi_\eps(\mx-\mz)=\varphi(\frac{x_1-z_1}{\eps^{\beta_1^j}},\dots,\frac{x_d-z_d}{\eps^{\beta_d^j}})$ in \eqref{j-1} instead of $\varphi$. Multiplying the obtained expression by $\eps$ and taking into account \eqref{sum1}, we conclude (only the principal symbols are in the sum below)

\begin{equation}
\label{j-2}
\sum\limits_{\malpha \in I'_j}(-1)^{|\malpha|}\sum\limits_{l=1}^m \int_{\R^d} \langle \rho_{j,\eps}(\mx-\my), a^{\malpha}_{jl}(\my) \mu_k(\my) \rangle \pa^\malpha \varphi_\eps({\mx-\mz}) d\mx+R_\eps=0,
\end{equation} where $R_\eps\to 0$ as $\eps\to 0$ (for clarifications, see the last term on the left-hand side in \eqref{i-3}). Next, we introduce the change of variables $\mx=(z_1+\eps^{\beta^j_1} w_1,\dots,z_d+\eps^{\beta_d^j}w_d)$ in the first term on the left-hand side of \eqref{j-2}. We get
\begin{equation}
\label{j-3}
\sum\limits_{\malpha \in I'_j}(-1)^{|\malpha|}\sum\limits_{l=1}^m \int_{\R^d} \langle \prod\limits^d_{k=1} \rho(\frac{z_k-y_k}{\eps^{\beta_k^j}}+w_k), a^{\malpha}_{jl}(\my) \mu_l(\my) \rangle \pa^\malpha \varphi(\mw) d\mw+R_\eps=0.
\end{equation} Now, we divide \eqref{j-3} by $|\mu_s|(E_\eps)$ and let $\eps\to 0$. Taking into account the uniform singularity assumptions as in \eqref{i-4}, we get 

\begin{equation}
\label{j-6}
\sum\limits_{\malpha \in I'_j}(-1)^{|\malpha|} \sum\limits_{l=1}^m  f_l(\mz) a^{\malpha}_{jl}(\mz) \int_{\R^{d}} {\mib \rho}(\mw)   \pa^\malpha \varphi(\mw) d\mw=0.
\end{equation} We now take
$$
\varphi(\mw)=\overline{T_{\frac{\psi(\pi_j(\mxi))}{|\xi_1|^{\beta^j_1}+\dots+|\xi_d|^{\beta^j_d}}} \rho (\mw)},
$$ where $T_{{m}}$ is the Fourier multiplier operator with the symbol ${m}$. After inserting this in \eqref{j-6} and applying the Plancherel theorem with respect to $\mw$, we obtain:
\begin{align}
\label{j-7}
\sum\limits_{l=1}^m  \sum\limits_{\malpha \in I'_j}(-1)^{|\malpha|} f_l(\mz) a^{\malpha}_{jl}(\mz)\int_{\R^{d}}   \frac{(2\pi \,i\mxi)^{\malpha}}{|\xi_1|^{\beta^j_1}+\dots+|\xi_d|^{\beta^j_d}}  \psi(\pi_j(\mxi)) |\hat{{\mib \rho}}(\mxi)|^2 d\mxi=0.&
\nonumber
\end{align} From here, due to arbitrariness of $\psi$, we get the claim.

{\bf Acknowledgement} The work is supported in part by the Croatian Science Foundation under Project WeConMApp/HRZZ-9780.


\end{document}